[1994/12/01]
\documentclass{article}
\usepackage{amsmath,amsthm}

\newtheorem{thm}{Theorem}[section]

\newtheorem{lem}[thm]{Lemma}

\theoremstyle{definition}

\newtheorem{conj}[thm]{Conjecture}

\newcommand{\BZ}{{\mathbf{Z}}}

\newcommand{\Lm}{{\Lambda}}

\newcommand{\floor}[1]{\left\lfloor#1\right\rfloor}

\theoremstyle{remark}

\title{A note on the paper by Bugeaud and Laurent "Minoration effective de la distance $p$-adique entre puissances
de nombres alg{\'e}briques"
\footnote{This paper is a revised version
of my master thesis\cite{Ymd2}.}}
\author{Tomohiro Yamada\footnote{Department of Mathematics, Kyoto University, Kyoto, 606-8502, Japan}\footnote{e-mail: \protect\normalfont\ttfamily{tyamada@math.kyoto-u.ac.jp}}}
  
\begin{document}

\maketitle

\begin{abstract}
We shall make a slight improvement to a result of $p$-adic logarithms,
which gives a nontrivial upper bound for the exponent of $p$ dividing
the Fermat quotient $x^{p-1}-1$.
\end{abstract}

\section{Introduction}

In 1939, Gel'fond\cite{Gel} established a result concerning upper bounds
for $p$-adic distance between two integral powers.  This result has been
refined by several papers such as Schinzel\cite{Sch}, Yu\cite{Yu1, Yu2, Yu3, Yu4, Yu5}, Bugeaud\cite{Bug}, and Bugeaud and Laurent\cite{BL}.
Our purpose is to improve a result in the last paper by noting that
constants $\gamma_j(j=1, 2)$ in \cite[Th{\'e}or{\`e}me 1]{BL} can be replaced by $1$,
which allows us to omit the condition $\log A_i\geq (\log p)/D$
in \cite{BL} by modifying some constants.

\begin{thm}\label{thm1}
Denote a $(p-1)$-th root of unity in $\BZ_p$ by $\zeta$.
Let $p$ be a prime and $\alpha_1, \alpha_2$ be integers not divisible by $p$.
Let $m_i$ be integers satisfying $\alpha_i\equiv\zeta ^{m_i}\pmod{p}$ for $i=1, 2$
and $g$ be an integer satisfying $\alpha_1^g\equiv\alpha_2^g\equiv 1\pmod{p}$.
Let $K\geq 3, L\geq 2$, $R_1, R_2, S_1, S_2$ be nonnegative integers.
Put $R=R_1+R_2-1, S=S_1+S_2-1, N=KL$.  Let $b_1, b_2$ be
positive integers with $(b_1, b_2, p)=1$ and denote
\begin{equation*}
b=\frac{(R-1)b_2+(S-1)b_1}{2}(\prod_{k=1}^{K-1}k!)^{-2/(K^2-K)}.
\end{equation*}
Suppose there exist congruence classes $c_1, c_2$ modulo $g$ such that
\begin{equation}\label{eqx1}
\begin{split}
&\#\{\alpha_1^r \alpha_2^s \mid 0\leq r<R_1, 0\leq s<S_1, m_1r+m_2s\equiv c_1\pmod{g}\}\\
&\geq L,
\end{split}
\end{equation}
\begin{equation}\label{eqx2}
\begin{split}
&\#\{b_2r+b_1s\mid 0\leq r<R_2, 0\leq s<S_2, m_1r+m_2s\equiv c_2\pmod{g}\}\\
&\geq (K-1)L.
\end{split}
\end{equation}
If we have
\begin{equation}\label{eqx3}
K(L-1)\log p>3\log N+(K-1)\log b+L((R-1)\log \alpha_1+(S-1)\log \alpha_2),
\end{equation}
then
\begin{equation}\label{eqx4}
v_p(\Lm)\leq KL-1.
\end{equation}
\end{thm}

A special case to which we can apply our version of Bugeaud-Laurent theorem is
a problem of Fermat quotient.  By a well-known theorem of Fermat, $x^{p-1} \equiv 1 \mod p$
for any prime $p$ and integer $x$ relatively prime to $p$.  However, it is unknown
whether there exist infinitely many prime $p$ such that $x^{p-1} \equiv 1 \mod p^2$.
It seems to be intersting and important to search for a nontrivial upper bound
for the exponent of $p$ dividing $x^{p-1}-1$.  This is equivalent to give a nontrivial
estimate for the $p$-adic logarithm $\log_p x^{p-1}$.  But already known results
for linear forms in $p$-adic logarithms do not give it.  As for results of
Bugeaud and Laurent\cite{BL}, the condition $\log A_i\geq (\log p)/D$ renders
the estimate trivial.  But now we can overcome this obstacle using Theorem \ref{thm1}.
Our result is as follows.

\begin{thm}\label{thm2}
If $p$ is a prime and $x, y$ are relatively prime integers, then
\begin{equation}\label{eq00}
v_p(x^{p-1}-1)\leq \floor{283(p-1)\frac{\log y}{\log p}\frac{\log xy}{\log p}}+4.
\end{equation}
If $q$ is odd prime, then we have
\begin{equation}\label{eq01}
v_p(q^{p-1}-1)\leq \floor{283(p-1)\frac{\log 2}{\log p}\frac{\log 2q}{\log p}}+4.
\end{equation}
Moreover, we have
\begin{equation}\label{eq02}
v_p(2^{p-1}-1)\leq \floor{283(p-1)\frac{\log 3}{\log p}\frac{\log 6}{\log p}}+4
\end{equation}
\end{thm}

Our argument is essentially the same as the argument of Bugeaud and Laurent\cite{BL}.
Indeed, all that we need is to make a very slight change in this paper.

Though this result is nontrivial, this seems to be far from best possible.
Ridout\cite{Rid} shows that there are only finitely many
rational integers $x$ such that $v_p(x^{p-1}-1)\geq (1+\epsilon)(\log x)/(\log p)$
for any fixed prime $p$ and positive $\epsilon$.  It is conjectured that
$v_p(x^{p-1}-1)\geq 3$ occurs only finitely many times for any fixed integer $x>1$.

The abc conjecture implies that for any $\epsilon>0$, the inequality
$v_p(x^{p-1}-1)\geq 1+(1+\epsilon p)(\log x)/(\log p)$ occurs only finitely
many times.  We can even conjecture:
\begin{conj}\label{conj1}
The inequality
\begin{equation}\label{eq03}
v_p(x^{p-1}-1)\leq 2+\frac{\log x+2\log\log x+\log\log p}{\log p}
\end{equation}
holds for any integer $x>1$ and prime $p$ except finitely many pairs $(x, p)$.
Furthermore, the inequality
\begin{equation}\label{eq04}
v_p(q^{p-1}-1)\leq 2+\frac{\log q+\log\log q+\log\log p}{\log p}
\end{equation}
holds for any primes $(q, p)$ except finitely many pairs $(q, p)$.
\end{conj}

We have a heuristic argument.  Since$(x+p)^{p-1}\equiv x^{p-1}-p\pmod{p^2}$,
we see that for any integer $x_0$ not divisible by $p$, the values
$(x^{p-1}-1)/p\pmod{p^{e-1}}(0\leq x<p^e-1, x\equiv x_0\pmod{p})$
take each congruent class exactly once.  Hence it is reasonable
to assume the probability of $x^{p-1}\equiv 1\pmod{p^e}$ is $p^{-e+1}$.
Let $e(x, p)$ be a function defined over nonnegative integers $x$ and
primes $p$.  If $\sum_{x, p}p^{-e(x, p)+1}$ converges, then we can expect
that $v_p(x^{p-1}-1)\geq e(x, p)$ has only finitely many solutions in $(x, p)$.

We can choose
\begin{equation}
e(x, p)=2+\frac{\log x+2\log\log x+\log\log p}{\log p}.
\end{equation}
Then we see that
\begin{equation}
\sum_{x, p}p^{-e(x, p)+1}\leq\sum_{x, p}(p\log p)^{-1}(x\log^2 x)^{-1}=\sum_{p}(p\log p)^{-1}\sum_{x}(x\log^2 x)^{-1}
\end{equation}
and therefore the sum converges.

One of our purposes of obtaining an upper bound for the exponent of $p$ dividing $x^{p-1}-1$
is an application for the study of problems involving the sum-of-divisors function.

Nagell\cite[Theorems 94, 95]{Nag} gives that $v_p(\sigma(q^c))\leq v_p(q^{p-1}-1)+v_p(c+1)$
for distinct primes $p, q$ with $q\ne 2$ and a positive integer $c$.  Now Theorem \ref{thm2}
immidiately gives the following theorem.
\begin{thm}\label{thm3}
If $q$ is an odd prime, then we have
\begin{equation}\label{eq05}
v_p(\sigma(q^c))\leq v_p(c+1)+\floor{283(p-1)\frac{\log 2}{\log p}\frac{\log 2q}{\log p}}+4.
\end{equation}
Moreover, we have
\begin{equation}\label{eq06}
v_p(\sigma(2^c))\leq v_p(c+1)+\floor{283(p-1)\frac{\log 3}{\log p}\frac{\log 6}{\log p}}+4.
\end{equation}
If we assume Conjecture \ref{conj1}, we have
\begin{equation}\label{eq07}
v_p(\sigma(q^c))\leq v_p(c+1)+2+\frac{\log q+\log\log q+\log\log p}{\log p}
\end{equation}
except only finitely many pairs $(p, q)$.
\end{thm}

We exhibit an application to the problem of perfect numbers.  If $N=\prod_{i=1}^{k}p_i^{e_i}$
is a perfect number with $p_1<\cdots<p_i$ distinct primes, then $e_i\leq (p_k-1)/2$ by a well-known
result of primitive prime factors.  Hence $N<(\prod_i p_i)^{(p_k-1)/2}$(for other finiteness results,
see, for example, \cite{Nie}, \cite{Pom}, \cite{Ymd1}).  We can improve this upper
bound using Theorem \ref{thm3}.
\begin{thm}\label{thm4}
If $\sigma(N)=\alpha N$ with $\alpha=n/d$ and $N=\prod_{i=1}^{k}p_i^{e_i}$ with $p_1<\cdots<p_i$ distinct primes, then
\begin{equation}\label{eq08}
N\leq d^C \prod_{i=1}^{k-1}p_i^{C(k-1)(p_k-1)/\log p_k}p_k^{C(k-1)(p_{k-1}-1)/\log p_{k-1}}
\end{equation}
for some absolute constant $C$.  Furthermore, if Conjecture \ref{conj1} is true, then
\begin{equation}\label{eq09}
N<\max\left\{e^{k^2}, d^{e/(e-2)}k^{-ke/(e-2)}(\prod_i p_i)^{kC'}\right\}
\end{equation}
for some absolute constant $C'$.
\end{thm}

We hope that our method will provide some systematical method to study
arithmetic functions involving divisors.

\section{Proof of Theorem \ref{thm1}}\label{prl}

We begin by improving \cite[Lemme 10]{BL}; we shall show that the term $g$
in the error terms can be omitted.
\begin{lem}
Let $K, L, R, S, g$ be integers $\geq 1$, $m_1, m_2, c$ rational integers with $(m_1, m_2, g)=1$.
Write $N=KL$ and $l_j=\floor{(j-1)/K} (j=1, \cdots, N)$.  If $(r_j, s_j) (j=1, \cdots, N)$
are $N$ pairs of integers satisfying
\begin{equation}\label{eqa1}
\begin{split}
0\leq r_v\leq R-1, 0\leq s_v\leq S-1,\\
m_1r_v+m_2s_v\equiv c\pmod{g}
\end{split}
\end{equation}
for $j=1, \cdots, N$, then
\begin{equation}\label{eqa2}
\begin{split}
M_1-G_1\leq \sum_{j=1}^{N}l_jr_j\leq M_1+G_1,\\
M_2-G_2\leq \sum_{j=1}^{N}l_js_j\leq M_2+G_2,
\end{split}
\end{equation}
where
\begin{equation}\label{eqa3}
\begin{split}
M_1=\frac{(L-1)(r_1+\cdots+r_N)}{2}, G_1=\frac{NL(R-1)}{4},\\
M_2=\frac{(L-1)(s_1+\cdots+s_N)}{2}, G_1=\frac{NL(S-1)}{4}.
\end{split}
\end{equation}
\end{lem}

\begin{proof}
We shall only show the inequality concerning $\sum_{j=1}^{N}l_jr_j$.  The other inequality
can be easily shown in a similar way.  Write $g'=(m_2, g)$ and $g''=g/g'$.

If $R\geq g'$, then we can proceed as in the original lemma and our lemma follows observing that
$R+g'-1\le 2(R-1)$.

If $R<g'$, then all $r_j$ must be equal to $c'$ in the original lemma.  Hence we have
\begin{equation*}
\sum_{j=1}^{N}l_jr_j=c'\sum_{j=1}^{N}l_j=c' N(L-1)/2=(r_1+\cdots +r_N)(L-1)/2=M_1.
\end{equation*}
This proves the lemma.
\end{proof}

In \cite[Lemme 11]{BL}, we can replace $\gamma_j(j=1, 2)$ by one.  This proves Theorem \ref{thm1}.

\section{Proof of Theorem \ref{thm2}}
We write $a_i=(\log A_i)/(\log p)$ and choose real constants $B, k, l$ satisfying $B\geq (\log b)/(\log p)$,
$k, l>0$, and
\begin{equation}\label{eq2a}
k\floor{lB+1}\floor{lB+2}-k\floor{lB+2}B\ge T_1+T_2+T_3,
\end{equation}
where
\begin{equation*}
\begin{split}
T_1=&2\floor{lB+2}^2\sqrt{k},\\
T_2=&2\floor{lB+2}^{3/2}(ga_1a_2)^{-1/2},\\
T_3=&\frac{3\log(ga_1a_2k\floor{lB+2}^2+\floor{lB+2})}{ga_1a_2\log p}.
\end{split}
\end{equation*}

We set $L, K, R_1, R_2, S_1, S_2$ as follows:
\begin{equation}
\begin{split}
L=\floor{lB}+2,& K=\floor{kgLa_1a_2}+1,\\
R_1=\floor{\sqrt{gLa_2/a_1}}+1,& S_1=\floor{\sqrt{gLa_1/a_2}}+1,\\
R_2=\floor{\sqrt{g(K-1)La_2/a_1}}+1,& S_2=\floor{\sqrt{g(K-1)La_1/a_2}}+1,\\
b=\frac{g(R+S-2)}{2}&(\prod_{k=1}^{K-1}k!)^{-2/(K^2-K)}.
\end{split}
\end{equation}

Proceeding as in \cite[Section 6]{BL}, we find that if the condition

\begin{equation}\label{eq2b}
\begin{split}
&\#\{b_2r+b_1s\mid 0\leq r<R_2, 0\leq s<S_2, m_1r+m_2s\equiv c\pmod{g}\}\\
=&\#\{(r, s)\mid 0\leq r<R_2, 0\leq s<S_2, m_1r+m_2s\equiv c\pmod{g}\}.
\end{split}
\end{equation}

holds, then we have $v(\Lm)<N$.

To prove Theorem \ref{thm2}, we apply this with $g=b_1=b_2=p-1$, $\alpha_1=q$, $\alpha_2=2, 3$
according to whether $p$ is odd or not, and each $m_i(i=1, 2)$ being an integer satisfying
$\alpha_i\equiv\zeta^{m_i}\pmod{p}$.  We may assume that $p>2^{283}$, since otherwise
the theorem follows from the trivial estimate.  We begin by confirming that the choice
$(k, l, B)=(11.32, 3, 1.027)$ satisfies (\ref{eq2a}).

Now $K=\floor{56.6ga_1a_2}+1>2^{100}$ and $\epsilon(K)$ in \cite[Lemme 13]{BL}
is smaller than $10^{-30}$.  Hence
\begin{equation}
\begin{split}
\log b&\leq\log g+\log(\frac{1}{a_1}+\frac{1}{a_2})+\frac{3}{2}-\log 2-\frac{1}{2}\log k+\epsilon(K)\\
&\leq\log p+\log\log p<(1.027)\log p,
\end{split}
\end{equation}
which assures that the choice $B=1.027$ satisfies the condition $B\geq (\log b)/(\log p)$.
Since $N=\floor{283ga_1a_2}+5$, we succeeded to prove (\ref{eq00}) under the condition (\ref{eq2b}).

If the condition (\ref{eq2b}) fails, then there exist rational integers $c_0, c_1$ such that
the congruence $m_1r+m_2s\equiv c_0\pmod{g}$ has two solutions $(r_1, s_1)$ and $(r_2, s_2)$ in integers
$0\leq r<R_2, 0\leq s<S_2$ satisfying $r_i-s_i=c_1$.  We have $r_2-r_1=s_2-s_1$ and therefore
\begin{equation}
(m_1-m_2)(r_1-r_2)\equiv m_1(r_1-r_2)+m_2(s_1-s_2)\equiv 0\pmod{g}.
\end{equation}
Let $m_0=m_1-m_2$ and $g_0$ be the residual order of $q\pmod{p}$.
We can easily see that $g_0$ divides $r_1-r_2$ since $\text{gcd}(m_1-m_2, g)=\text{gcd}(m_0, g)=g/g_0$.  Hence $R_2>\max\{r_1, r_2\}\geq g_0$.
Now (\ref{eq00}) follows from the trivial estimate
\begin{equation}
v(Lm)\leq g_0\frac{\log q}{\log p}\leq (R_2-1)\frac{\log q}{\log p}\leq gLa_2{\log q}{\log p}=5ga_1a_2.
\end{equation}

(\ref{eq01}) and (\ref{eq02}) immidiately follow from (\ref{eq00}) by
taking $(x, y)$ as $(2, q)$ and $(3, 2)$ respectively.  This completes
the proof.

\section{Proof of Theorem \ref{thm4}}
In this section, we denote by $C_1, C_2, \cdots $ absolute constants.

By Theorem \ref{thm3}, we have
\begin{equation}
v_{p_j}(\sigma(p_i^{e_i}))\leq v_{p_j}(e_i+1)+283(p_j-1)\frac{\log 2}{\log p_j}\frac{\log 2p_i}{\log p_j}+4
\end{equation}
if $p_i$ is odd.  Moreover, we have
\begin{equation}
v_{p_j}(\sigma(2^{e_i}))\leq v_{p_j}(e_i+1)+283(p_j-1)\frac{\log 3}{\log p_j}\frac{\log 6}{\log p_j}+4.
\end{equation}
Thus we obtain
\begin{equation}
\begin{split}
v_{p_j}(\alpha N)&\leq v_{p_j}(\prod_{i\neq j}(e_i+1))\\
&+283\frac{p_j-1}{(\log p_j)^2}(\sum_{p_i=2}\log 3\log 6+\log 2\sum_{i\neq j, p_i>2}\log 2p_i)\\
&+4(k-1).
\end{split}
\end{equation}
Noting that $N$ divides $d\alpha N$ and $N$ is composed of $p_i(1\leq i\leq k)$, we see
\begin{equation}
\begin{split}
N&\leq d\prod_{i=1}^{k}(e_i+1)\prod_{p_i=2}6^{4(k-1)+(283\log 3)\sum_{p_j>2}(p_j-1)/\log p_j}\\
&\times\prod_{p_i>2}(2p_i)^{4(k-1)+(283\log 2)\sum_{j\neq i}(p_j-1)/\log p_j}.
\end{split}\end{equation}
Since $\prod_{i=1}^{k}(e_i+1)=d(N)=N^{o(1)}$, we have
\begin{equation}
N\leq d^{C_1} \prod_{i=1}^{k-1}p_i^{C_1(k-1)(p_k-1)/\log p_k}p_k^{C_1(k-1)(p_{k-1}-1)/\log p_{k-1}}.
\end{equation}
This proves (\ref{eq08}).

We assume Conjecture \ref{conj1}.  Theorem \ref{thm3} gives
\begin{equation}
v_{p_j}(\sigma(p_i^{e_i}))\leq v_{p_j}(e_i+1)+\max\{2+\frac{\log p_i+\log\log p_i+\log\log p_j}{\log p_j}, C_2\}
\end{equation}
for any $i, j$.  A similar argument to the first case, we have
\begin{equation}\label{eq31}
\begin{split}
N&\leq d\prod_{i=1}^{k}(e_i+1)\prod_{1\leq i, j\leq k, i\neq j} \max\{p_i(\log p_i)p_j^2(\log p_j), p_j^{C_2}\}\\
&\leq d\prod_{i=1}^{k}(e_i+1)\prod_{1\leq i, j\leq k, i\neq j}p_i(\log p_i)p_j^{C_2}(\log p_j)\\
&\leq d\prod_{i=1}^{k}(e_i+1)p_i^{(C_2+1)(k-1)}(\log p_i)^{2(k-1)}\\
&\leq d\prod_{i=1}^{k}(e_i+1)p_i^{C_3(k-1)}.
\end{split}
\end{equation}
Let $E=\sum(e_i+1)\log p_i$.  Then, since
\begin{equation}
\prod_{i=1}^{k}p_i^{e_i-C_3(k-1)}\leq d\prod_{i=1}^{k}(e_i+1)
\end{equation}
by (\ref{eq31}), we have
\begin{equation}
\prod_{i=1}^{k}(\log p_i)p_i^{e_i+1-kC_4}\leq d\prod_{i=1}^{k}(e_i+1)\log p_i\leq d(E/k)^k.
\end{equation}
Taking the logarithms of both sides, we have
\begin{equation}
E\leq k\log E-k\log k+\log d+kC_4\sum_i\log p_i.
\end{equation}
We observe that $x\geq k^2$ implies $k(\log x)/x\leq 2(\log k)/k\leq 2/e$.  Hence
\begin{equation}
E\leq \max\left\{k^2, \frac{e}{e-2}(\log d-k\log k+kC_4\sum_i\log p_i)\right\}.
\end{equation}
Since $N=\prod_{i} p_i^{e_i}<e^E$, we have
\begin{equation}
N<\max\left\{e^{k^2}, d^{e/(e-2)}k^{-ke/(e-2)}(\prod_i p_i)^{kC_5}\right\}.
\end{equation}
This completes the proof.

{}
\end{document}